\documentclass [12 pt]{article}
\usepackage{amsfonts}
\def\lr{\left(}
\def\rr{\right)}
\def\a{\alpha}
\def\lc{\left\{}
\def\rc{\right\}}
\def\qed{\vbox{\hrule\hbox{\vrule\kern3pt\vbox{\kern6pt}\kern3pt\vrule}\hrule}}
\newtheorem{thm}{Theorem}
\newtheorem{lem}[thm]{Lemma}

\title{$t$-Covering Arrays Generated by a Tiling Probability Model}
\author{Michael S.~Donders \\
Department of Mathematics and Computer Science\\
McDaniel College
\and
Anant P.~Godbole\\
Department of Mathematics and Statistics\\
East Tennessee State University}
\begin{document}
\def\nl{\newline}
\def\N{ n(m,t,\alpha)}
\maketitle
\begin{abstract}  A $t-\a$ covering array is an $m\times n$ matrix, with entries from an alphabet of size $\alpha$, such that for any choice  of $t$ rows, and any ordered string of $t$ letters of the alphabet, there exists a column such that the ``values" of the rows in that column match those of the string of letters. We use the Lov\'asz Local Lemma in conjunction with a new tiling-based probability model to improve the upper bound on the smallest number of columns $N=N(m,t,\alpha)$ of a $t-\a$ covering array.
\end{abstract}

\section
{Introduction}  Consider an $m\times n$ matrix with entries from the ``alphabet"  $A=\{1, 2, \ldots, \alpha \}$. Let the $(i,j)^{\rm th}$ entry be represented by $r_{i,j}$. We say that this matrix is a {\em $t-\alpha$-covering matrix} or a  {\em $t-\alpha$-covering array} if given any $t$ rows, $p_1, p_2, \ldots ,p_t$ of the matrix, and any vector $\langle v_1, v_2, \dots ,v_t\rangle$, with $v_i \in A$, there exists a column $q$ such that
$$\langle v_1 ,v_2, \ldots , v_t\rangle = \langle r_{p_1,q},r_{p_2,q}, \ldots ,r_{p_t,q}\rangle.$$  Extensive surveys of covering arrays may be found in the papers of Sloane \cite{sloane} and Colbourn \cite{cc}.  
Given $t, m$ and the alphabet size $\vert A\vert$, we wish to find the minimum number of columns, $n$, such that there exists an $m\times n$ matrix that is $t$-covering. We will define $N=N(m,t,\alpha)$ as the smallest positive integer $n$ such that there exists a covering array of dimensions $m\times n$.  At the Coimbra Zero-One Matrix Conference, the second author talked about the need to introduce new probability models to improve upper bounds on $N(m,t,\alpha)$ and the corresponding numbers for {\it partial covering arrays} \cite{patty}.  In this paper we propose a specific way of doing so, once again using the Lov\'asz local lemma as an auxiliary tool.
\begin{lem}
 The Lov\'asz Local Lemma (\cite{alon}):  
Let $C_1,C_2, \dots,\ C_K$ be the events in arbitrary probability space. Suppose that each event $C_i$ is mutually independent of a set of all the other events $C_k$ but at most $d$, and that $P(C_i) \le p$ for all $1\le i \le K$. If $ep(d+1) \le 1$ then $P(\bigcap^K_{k=1} C'_k) > 0.$
\end{lem}

 Let $R$ be the index set of all sets of $t$ rows; $\vert R\vert = {m\choose t}$. For $r \in R$, let $C_r$ be the event that the $r^{\rm th}$ row set does not contain some vector $\langle v_1 ,v_2, \ldots , v_t\rangle$ in any of its columns.  We wish to prove that $P(\bigcap_{r\in R}C'_r)>0$ if $n\ge N_0$, proving that $N(m,t,\alpha)\le N_0$.   Now in \cite{gss} a general upper bound was provided on the size of covering arrays; this was
\begin{equation}N(m,t,\alpha)\le N_0:= (t-1)\frac{\log_2(m)}{\log_2\frac{\alpha^t}{\alpha^t-1}}\{1 + o(1)\}.\end{equation}
The proof used an elementary probability model that consisted of placing one letter of the alphabet independently in each of the $mn$ positions with probability $\frac{1}{\alpha}$, i.e. by letting $P(r_{i,j} = x) = \frac{1}{\alpha}\ \forall x\in A$.  In the same paper, a special probability model was used, but only for the case $\alpha=2, t=3$.  Here the authors of \cite{gss}, following the approach used in the doctoral thesis of Roux (see, e.g. \cite{sloane}), used a probability model that independently places an equal number of zeros and ones in the rows of the matrix (the so-called ``fixed weight rows" model.)  Unfortunately this method becomes quite intractable in general, and it is our intent in this paper to explore a probability model that is, in some sense, intermediate between the general technique in \cite {gss} and the special method used there for $\alpha=2,t=3$:
Specifically, we seek to improve the general bound (1) using the method of placing {\it consecutive and equally weighted tiles} along the rows. We use tiles of dimension $1\times k\alpha$, such that there are exactly $k$ $x$'s in each tile for each $x\in A$.  By way of comparison, the general method used $1\times 1$ tiles that led to a loss of control over the numbers of letters of each type in any row, while Roux's method used a {\it single long tile}  in each row, i.e., corresponded to $k=n/2$ ($n$ even).  

We consider  two cases, when (i) $k=1$ which yields an elementary equation relating $N(m,t,\alpha)$ and the variables $m$, $t$ and $\alpha$, and when (ii) $k > 1$, which yields better bounds as $k$ increases, but which generates increasingly more complicated solutions.  

\noindent (i) We start with the case $k = 1$, and fill in our matrix using tiles that contain one randomly placed copy of each letter of the alphabet, assuming that $\alpha\vert n$. Note that there are a total of $\alpha^t$ possible vectors, and by the symmetry of our construction, all are equally likely to occur in the selected rows. Thus $P(C_r) \le \lambda \alpha^t$ where $\lambda$ is the probability that a specific vector $z^* = \langle z_1 ,z_2,\ldots, z_t\rangle$ is missing in the set $r$ of selected rows. Select an arbitrary set of $t$ rows in the  matrix. Consider the columns in any vertically aligned set of tiles.  For each $z_i$, there is exactly one value in any tile equal to $z_i$, and $\alpha$ places it can be; moreover $z^*$ cannot occur in more than one column of the vertically stacked tiles in the selected rows. Therefore, the probability that $z^*$ is somewhere in these tiles is $\alpha \cdot\lr\frac{1}{\alpha}\rr^t = \lr\frac{1}{\alpha}\rr^{t-1} $. Since there are $\frac{n}{\alpha}$ tiles in any row of the $m\times n$ matrix, and the composition of these is determined independently, we have 

$$\lambda = \lr 1- \lr\frac{1}{\alpha}\rr^{t-1}\rr ^{\frac{n}{\alpha}} = \lr\lr1-{\lr\frac {1}{\alpha}\rr}^{t-1}\rr ^{\frac{1}{\alpha}}\rr ^n,$$
and thus, 
\begin{equation}P(C_r) \le \alpha^t\lr\lr{1 - \lr\frac {1}{\alpha}\rr}^{t-1}\rr ^{\frac{1}{\alpha}}\rr ^n.\end{equation}

\noindent We can improve this bound slightly by using a technique found in \cite{patty}, where the vectors $z_i = \langle i,i, \ldots , i\rangle; 1\le i\le \alpha$ can be achieved for all sets $ r$ by including columns consisting of all $i$'s. There are $\alpha$ of these vectors; thus this reduces the number of $z^*$'s from $\alpha^t$ to $\alpha^t - \alpha$. We can ignore these vectors in our calculation of $P(C_r)$ so long as we remember to add $\alpha$ columns to our value $N(m,t,\alpha)$. So (2) may be improved as follows: 
\begin{equation}P(C_r)\le (\alpha^t-\alpha)\lr\lr{1 - \lr\frac {1}{\alpha}\rr}^{t-1}\rr ^{\frac{1}{\alpha}}\rr ^n.
\end{equation}
Our next step is to calculate $d$. For any set $r$ of rows, there will be a dependency only on sets $r_0 \in R$ such that $r \cap r_0 \neq \emptyset$. We will bound the number of such $r_0$'s by choosing one row from $r$, and then choosing an arbitrary $t-1$ rows from the $m-1$ other rows in the matrix.  Thus $d \le t{m-1 \choose t-1}$, so $d+1 \le \frac{t m^{t-1}}{(t-1)!}$. Substituting this into the Lov\'asz local lemma we get $$ep(d+1) = \frac{e t m^{t-1}}{(t-1)!}\lr\alpha^t-\alpha\rr\lr\lr{1-\lr\frac {1}{\alpha}\rr}^{t-1}\rr ^{\frac{1}{\alpha}}\rr ^n \le 1$$
if
$$n \ge \frac{(t-1)\log_2(m)}{\log_2\lr\lr\frac{\alpha^{t-1}}{\alpha^{t-1}-1}\rr^{\frac{1}{\alpha}} \rr}\lc 1 + \frac{\log_2\lr\alpha^t - \alpha\rr}{(t-1)\log_2(m)} + \frac{\log_2(et)}{(t-1)\log_2(m)} -  \frac{\log_2\lr(t-1)!\rr}{(t-1)\log_2(m)} \rc,$$
i.e., if
$$ n \ge \frac{\alpha(t-1)\log_2(m)}{\log_2\lr\frac{\alpha^{t-1}}{\alpha^{t-1}-1}\rr}\{1 + o(1)\}\ \ \ \ \ \  m \to \infty.$$
It follows that
\begin{equation} N(m,t,\alpha) \le \frac{\alpha(t-1)\log_2(m)}{\log_2\lr\frac{\alpha^{t-1}}{\alpha^{t-1}-1}\rr}\{1 + o(1)\}.\end{equation}
since adding back, into (4), the $\alpha$ columns we removed earlier only changes the $o(1)$ term.  Notice that the above process gives us  both a precise and an asymptotic bound for $\N $. Note too that (4) gives an improvement over the previous best bound (1) due to the fact that 
\[\lr1-\frac{1}{\alpha^{t-1}}\rr^{\frac{1}{\alpha}}\le1-\frac{1}{\alpha^t}.\]

\noindent (ii) We now consider the case $k > 1$; recall that the size of our tiles is $1\times k\alpha$.
First note that the size of the tile does not change $d$, and thus $d+1 \le \frac{t m^{t-1}}{(t-1)!}$ as before. We next reconsider $P(C_r)$, and compute it using inclusion exclusion. Let 
\begin{eqnarray}\gamma_k &=&\frac{ \sum_{i=1}^k(-1)^{(i+1)}{\alpha k \choose i}{{\alpha k - i} \choose {k-i,k,\ldots,k}}^t}{{{\alpha k} \choose {k,k,\ldots,k}}^t}\nonumber\\
&=&\frac{ \sum_{i=1}^k(-1)^{(i+1)}{\alpha k \choose i}{\alpha k - i \choose k-i}^t}{{\alpha k \choose k}^t}\end{eqnarray}
be the probability that a given vector $z^*$ is in a given vertical array of $t$ tiles. This yields $\lambda=\lambda_k = (1-\gamma_k)^\frac{n}{k\alpha}$ and hence
$$ P(C_r)\le (\alpha^t-\alpha)(1-\gamma_k)^\frac{n}{k\alpha},$$
so that the Lov\'asz local lemma yields $P(\cap C'_r)>0$ if,
$$(\alpha^t-\alpha)(1-\gamma_k)^\frac{n}{k\alpha} \frac{et m^{t-1}}{(t-1)!}\le 1$$
i.e., if ,
$$n \ge \frac{(t-1)\log_2(m)}{\log_2\lr\lr\frac{1}{1-\gamma_k}\rr^{\frac{1}{k\alpha}} \rr}\lc 1 + \frac{\log_2\lr\alpha^t - \alpha\rr}{(t-1)\log_2(m)} + \frac{\log_2(et)}{(t-1)\log_2(m)} -  \frac{\log_2\lr(t-1)!\rr}{(t-1)\log_2(m)} \rc,$$
or
$$n \ge \frac{k\alpha (t-1)\log_2(m)}{\log_2\lr\frac{1}{1-\gamma_k}\rr}\{1 + o(1)\} \ \ \ \ \  \ \ m\rightarrow \infty.$$  It follows that
$$N(m,t,\alpha) \le \frac{k\alpha (t-1)\log_2(m)}{\log_2\lr\frac{1}{1-\gamma_k}\rr}\{1 + o(1)\}.\eqno(6)$$

\noindent{\bf Comments}  It is clear that as we increase $k$ from 1 to $\frac{n}{\alpha}$, the bound on $N(m,t,\alpha)$ becomes better and better, while the equation to solve for it becomes more and more convoluted. 
Take for example, the case when $t=3$, $\alpha = 2$. The previous best known bound (1) for a general $N(m,t,\alpha)$ yields the solution $N(m,3,2) \le 10.38 \log_2(m)\{1+o(1)\}$, while the best known solution for this specific case (Roux [4]) yields $N(m,3,2) \le 7.56 \log_2(m)\{1 + o(1)\}$, a result obtained by equally weighing all the rows to have the same number of 1's and 0's. The solution obtained via tiling yields $N(m,3,2) \le 9.64\log_2(m)\{1 + o(1)\}$ when $k=1$.  With this we can see that even the simplest case of the tiling solution, $k=1$, offers a fairly significant improvement in the bounds, while more complex solution will provide the better bounds for the size of a covering array.  A few values of $N(m,t,\alpha)$ as given by (4) and (6) may be found in the following table; $k=0$ refers to the bound in (1):\vfill\eject
$$\vbox{\halign{
\hfil#\hfil&\qquad
\hfil#\hfil&\qquad
\hfil#\hfil&\qquad
\hfil#\hfil\cr
$\alpha$&	$t$&	$k$&	$N(m,t,\alpha)/\log_2m$	\cr
2&	3&	0&	10.38\cr
2&	3&	1&	9.64\cr
2&	3&	2&	8.68\cr
2&	3&	3&	8.31\cr
2&	4&	0&	32.22\cr
2&	4&	1&	31.15\cr
2&	4&	2&	29.55\cr
2&	4&	3&	28.85\cr
3&	3&	0&	36.73\cr
3&	3&	1&	35.31\cr
3&	3&	3&	33.28\cr
3&	3&	5&	32.79\cr
3&	4&	0&	167.39\cr
3&	4&	1&	165.3\cr
3&	4&	3&	161.57\cr
3&	4&	5&	160.47\cr
4&	3&	0&	88.03\cr
4&	3&	2&	83.97\cr
4&	3&	4&	82.72\cr
4&	3&	6&	82.27\cr
4&	4&	0&	531.3\cr
4&	4&	2&	524.75\cr
4&	4&	4&	521.98\cr
4&	4&	6&	520.90\cr
5&	4&	0&	1298.61\cr
5&	4&	2&	1290.12\cr
5&	4&	4&	1286.46\cr
5&	4&	6&	1285.01\cr
5&	5&	0&	8662.95\cr
5&	5&	2&	8651.13\cr
5&	5&	4&	8644.67\cr
5&	5&	6&	8641.86\cr
}}$$

\section{Open Problems}  Perhaps the overarching open problem is that of using alternative probability models in order to tease out better and better bounds on the size of minimal covering arrays.  Markov models and others involving global dependence are one option.  A method more relevant to the central problem addressed at the Coimbra conference, would, however, be to work with zero-one or alphabet based matrices with fixed row and column totals (in this paper we fix just the row totals!).  Last but not least, can we let $k$ go to infinity (at a relatively slow rate) and analyze the sum in (5)?  Can we conduct the analysis with $k=n/\alpha$?
\section {Acknowledgment}  The research of both authors was supported by NSF Grant 1004624.

\end{document}